\documentclass[11pt,a4paper]{amsart}
\usepackage{amssymb}
\usepackage[T1]{fontenc}
\usepackage{Alegreya}
\usepackage[libertine]{newtxmath}
\usepackage[headings]{fullpage}
\usepackage{paralist}
\usepackage{tikz}

\usepackage{color}

\usepackage[colorlinks=true,linkcolor=blue, citecolor=blue, urlcolor=blue,%
pagebackref=true]{hyperref}
\makeatletter

\@addtoreset{equation}{section}
\def\theequation{\thesection.\@arabic \c@equation}

\def\theenumi{\@alph\c@enumi}

\makeatother
\theoremstyle{plain}
\newtheorem{theorem}[equation]{Theorem}

\newtheorem{corollary}[equation]{Corollary}
\newtheorem{proposition}[equation]{Proposition}

\theoremstyle{definition}

\newtheorem{remark}[equation]{Remark}

\newtheorem{example}[equation]{Example}

\newtheorem{definition}[equation]{Definition}

\newtheorem{notation}[equation]{Notation}

\newtheorem{discussion}[equation]{Discussion}
\newenvironment{discussionbox}[1][]{%
    \begin{discussion}[#1]\pushQED{\qed}}{\popQED \end{discussion}}

\newtheorem{observation}[equation]{Observation}

\newtheorem{construction}[equation]{Construction}

\newtheorem{setup}[equation]{Setup}
    {\setcounter{step}{0}}{} \newcounter{step}

\newcommand{\calM}{\mathcal M}

\newcommand{\naturals}{\mathbb{N}}
\newcommand{\ints}{\mathbb{Z}}

\newcommand{\reals}{\mathbb{R}}

\def\to{\longrightarrow}
\newcommand{\define}[1]{\emph{#1}}
\newcommand{\minus}{\ensuremath{\smallsetminus}}
\DeclareMathOperator{\JI}{JI}
\DeclareMathOperator{\Des}{Des}
\DeclareMathOperator{\Pow}{Pow}
\DeclareMathOperator{\maxChains}{\calM}

\newcommand{\arrowIn}{
	\tikz \draw[-stealth] (-1pt,0) -- (1pt,0);
}

\begin{document}
\title{The $h$-polynomial and the rook polynomial of some polyominoes}

\author{Manoj Kummini}
\address{Chennai Mathematical Institute, Siruseri, Tamilnadu 603103. India}
\email{mkummini@cmi.ac.in}

\author{Dharm Veer}
\address{Chennai Mathematical Institute, Siruseri, Tamilnadu 603103. India}
\email{dharm@cmi.ac.in}

\thanks{MK was partly supported by the grant CRG/2018/001592
from Science and Engineering Research Board, India and
by an Infosys Foundation fellowship.
DV was partly supported by an Infosys Foundation fellowship.}
\begin{abstract}
Let $X$ be
a convex polyomino such that its vertex set is a sublattice 
of $\ \naturals^2$. Let $\Bbbk[X]$ be the toric ring (over a field $\Bbbk$)
associated to $X$ in the sense of Qureshi, \emph{J. Algebra}, 2012.
Write the Hilbert series of $\Bbbk[X]$ as
$(1 + h_1 t + h_2 t^2 + \cdots )/(1-t)^{\dim(\Bbbk[X])}$. For $k \in
\naturals$, let $r_k$ be the number of configurations in $X$ with $k$
pairwise non-attacking rooks. 
We show that $h_2 < r_2$ if $X$ is not a thin polyomino. 
This partially confirms a conjectured characterization of thin
polyominoes by Rinaldo and Romeo, \emph{J. Algebraic Combin.}, 2021.

\end{abstract}

\maketitle

\section{Introduction}
\label{section:intro}

A polyomino is a finite union of unit squares with vertices at lattice
points in the plane that is connected and has not finite 
cut-set~\cite[4.7.18]{StanEC1}.
(Definitions are in Section~\ref{section:prelims}.)
A.~A.~Qureshi~\cite{QureshiPolyominoes2012} 
associated a finitely generated graded algebra $\Bbbk[X]$
(over a field $\Bbbk$) to polyomino $X$.
For $k \in \naturals$, a $k$-rook configuration in $X$ is an arrangement of
$k$ rooks in pairwise non-attacking positions. The rook polynomial $r(t)$
of $X$ is $\sum_{k \in \naturals } r_k t^k$ where $r_k$ is the number of
$k$-rook configurations in $X$. The $h$-polynomial of $\Bbbk[X ]$ is the 
(unique) polynomial $h(t) \in \ints[t]$ such that the Hilbert series of 
$\Bbbk[X ]$ is $h(t)/(1-t)^d$ where $d = \dim \Bbbk[X ]$. A polyomino is
thin if it does contain a $2 \times 2$ square of four unit squares (such
as the one shown in Figure~\ref{figure:twoRooks}).

G.~Rinaldo and
F.~Romeo~\cite[Theorem~1.1]{RinaldoRomeoHilbSeriesThinPolyominoes2021}
showed that if $X$ is a simple thin polyomino, then $h(t) = r(t)$ and 
conjectured~\cite[Conjecture~4.5]{RinaldoRomeoHilbSeriesThinPolyominoes2021} 
that this property characterises thin
polyominoes. In this paper, we prove this conjecture in the following case:
\begin{theorem}
\label{theorem:RRconj}
Let $X$ be a convex polyomino such that its vertex set $V(X )$ is a sublattice 
of $\ \naturals^2$. 
Let $h(t) = 1 + h_1 t + h_2 t^2 + \cdots$ be the $h$-polynomial of 
$\Bbbk[X]$ and $r(t) = 1 + r_1 t + r_2 t^2 + \cdots$ be the rook polynomial
of $X$.
If $X$ is not thin, then $h_2 < r_2$.
In particular $h(t) \neq r(t )$.
\end{theorem}

Its proof proceeds as follows: 
we first observe that $\Bbbk[X]$ is the Hibi ring of the distributive
lattice $V(X)$ and
that the Hilbert series of the Hibi ring of a distributive lattice and that
of the Stanley-Reisner ring of its order complex are the same. We then use
the results of~\cite{BjornerGarsiaStanley1982} relate the $h$-polynomial to
descents in maximal chains of $V(X)$, and find an injective map from
the set of maximal chains of $V(X)$ to the rook configurations in $X$, to
conclude that $h_k \leq r_k$ in general. We
then show that if $X$ is not thin, this map is not surjective to show that
$h_2 < r_2$. In Corollary~\ref{corollary:Lconvex} 
we extend our result to $L$-convex polyominoes.

Section~\ref{section:prelims} contains the definitions and preliminaries.
Proof of the theorem is given in Section~\ref{section:proof}.

\section*{Acknowledgements}
The computer algebra systems Macaulay2~\cite{M2} and
SageMath~\cite{sagemath} provided valuable assistance in studying examples.

\section{Preliminaries}
\label{section:prelims}

\begin{definition}
A \define {cell} in $\reals^2$ is a set of the form $\{(x, y) \in \reals^2 \mid
a \leq x \leq a+1, b\leq y \leq b+1\}$ where $(a,b ) \in \ints^2$. 
We identify the cells of $X$ by their top-right corners: For $v \in
\ints^2$, $C(v )$ is the cell whose top-right corner is $v$.
A \define{polyomino} $X$ is a finite union of cells that is connected and
has no finite
cut-set (i.e., removing finite sets from $X$ leaves $X$ 
connected)~\cite[4.7.18]{StanEC1}.
We say that a polyomino $X$ is \define{horizontally convex} 
if for every line segment $\ell$ parallel to the $x$-axis with end-points
in $X$, $\ell \subseteq X$. Similarly we define 
\define{vertically convex} polyominoes.
We say that a polyomino $X$ is \define{convex}
if it is horizontally convex and vertically convex.
The set of cells of $X$ is denoted by $C(X)$.
The \define{vertex set} $V(X)$ of $X$ is $X \cap \ints^2$.
By \define the {left-boundary vertices} of $X$, we mean the 
elements of $\ints^2 \cap \partial X$ that are
top-left vertices of the cells of X; 
the \define {bottom-boundary vertices} of $X$ are the 
elements of $\ints^2 \cap \partial X$ that are
bottom-right vertices of the cells of X; 
\end{definition}

Qureshi~\cite{QureshiPolyominoes2012} associated a toric ring to a
polyomino.

\begin{definition}
Let $X$ be a convex polyomino. Let $R = \Bbbk[\{x_v \mid v \in V(X)\}]$ be a
polynomial ring. An \define{interval} in $X$ is a subset of $X$ of the form
$[a,b ] := \{c \in V(X ) \mid a \leq c \leq b \}$ where $a \leq b \in V(X)$ 
and $\leq$ is the
partial order on $\reals^2$ given by componentwise comparison: $a = (a_1,
a_2) \leq  b = (b_1, b_2 )$ if $a_1 \leq b_1$ and $a_2 \leq b_2$. Let $I_X$
be the $R$-ideal generated by the binomials of the form 
$x_a x_b - x_c x_d$ where $a \leq b \in V(X )$ and $c, d \in V(X )$ are the
other two corners of the interval $[a,b ]$.
Let $\Bbbk[X ] = R/I_X$.
\end{definition}

\begin{setup}
\label{setup:convexpolyomino}
Let $X$ be a convex polyomino such that $V(X )$ is a sublattice of
$\naturals^2$. 
Let $\JI(X )$ be the poset of join-irreducible elements of $V(X )$.
After a suitable translation, if necessary, we assume that
$(0,0)$ and $(m,n )$ are the elements $\hat 0$ and $\hat 1$ of $V(X )$.
Hence $|\JI(X ) | = m+n$. 
\end{setup}

\begin{definition}
Let $L$ be a finite distributive lattice.
Let $R = \Bbbk[\{x_a \mid a \in L \}]$. The
\define{Hibi ideal}~\cite{HibiDistrLatticesASLs1987} 
$I_L$ of $L$ is the $R$-ideal generated by the
binomials of the form $x_ax_b - x_cx_d$ where $a,b \in L$ 
and $c$ and $d$ are the join and the meet of $a$ and $b$.
The \define{Hibi ring} of $L$ is $\Bbbk[L] := R/I_L$.
\end{definition}

\begin{definition}
Let $R$ be a standard graded $\Bbbk$-algebra. The \define{$h$-polynomial}
of $R$ is the polynomial $h(t)$ such that the Hilbert series of $R$ is 
$h(t)/(1-t)^d$ where $d =\dim R$.
\end{definition}

\begin{remark}
\label{remark:hPolyHibiRingSRRing}
When $X$ is as in Setup~\ref{setup:convexpolyomino},
the polyomino ring $\Bbbk[X]$ is the Hibi ring $\Bbbk[V(X)]$.
Hence we are interested in the $h$-polynomial of the Hibi ring of a
distributive lattice. 
Let $L$ be a distributive lattice.
The order complex $\Delta(L)$ is the simplicial complex whose faces are the
chains of $L$. The \define{Stanley-Reisner ring} $\Bbbk[\Delta(L)]$
of $\Delta(L)$ is the
quotient of $\Bbbk[\{x_a \mid a \in L \} ]$ by the ideal generated by
$\{x_ax_b \mid a, b \;\text{incomparable}\}$.
There is a flat deformation from $\Bbbk[L]$ to $\Bbbk[\Delta(L)]$; see,
e.g.,~\cite[Section~7.1]{BrHe:CM}, after noting that Hibi rings are ASLs.
Hence the $h$-polynomials of $\Bbbk[X]$ and of $\Bbbk[\Delta(V(X))]$
are the same. We use the results of~\cite{BjornerGarsiaStanley1982}
to relate the $h$-polynomial of $\Delta(L)$ to the descents in the maximal
chains of $L$.
\end{remark}

\begin{discussionbox}
\label{discussionbox:descent}
We follow the discussion of~\cite[Section~1]{BjornerGarsiaStanley1982}.
Let $\omega : \JI(X ) \to \{1, \ldots, m+n \}$ be 
a (fixed) order-preserving map.
Let $\maxChains(X)$ be the set of maximal chains of $V(X )$.
Let $\mu \in \maxChains(X)$. We first write $\mu$ as a chain of order
ideals of $\JI(X)$:
$\hat 0 = I_0 \subsetneq I_1 \subsetneq \cdots \subsetneq I_{m+n} = \hat 1$.
Then $|I_i \minus I_{i-1}| = \{p_i\}$ for some $p_i \in \JI(X )$.
Define $\omega(\mu ) = (\omega (p_1 ), \ldots, \omega(p_{m+n }))$.
For $1 \leq i \leq m+n-1$,
we say that $i$ is a \define{descent} of $\mu$
if $\omega( p_i) > \omega (p_{i+1})$.
The \define{descent set} $\Des(\mu)$ of $\mu$ is 
$\{i \mid 1 \leq i \leq m+n-1, \;i \;\text{is a descent of}\; \mu \}$.
For $k \in \naturals$, 
define $\maxChains_k(X) = \{\mu \in \maxChains(X) : |\Des(\mu) |=k\}$.

We now think of $\mu$ as a lattice path from $(0,0 )$ to $(m,n)$ consisting of
horizontal and vertical edges. Label the vertices of $\mu$ as $(0,0) = \mu_0,
\mu_1, \ldots, \mu_{m+n } = (m,n )$, with $\mu_i - \mu_{i-1 }$ a unit
vector (when we think of these as elements of $\reals^2$)
pointing to the right or upwards. Then, if $i \in \Des(\mu )$, then the
direction of $\mu$ changes at $\mu_i$, i.e, the vectors $\mu_i-\mu_{i-1 }$
and $\mu_{i+1 }-\mu_i$ are perpendicular to each other. Hence $\mu_{i-1 }$
and $\mu_{i+1 }$ are the bottom-left and top-right vertices of a cell 
(the cell $C(\mu_{i+1})$ in our notation)
of $X$. Thus we get a function
\begin{equation}
\label{equation:MXtoPowCX}
\psi: \maxChains(X ) \to \Pow(C(X)), \qquad 
\mu \mapsto \{C(\mu_{i+1}) \in C(X) \mid i \in \Des(\mu) \}.
\qedhere
\end{equation}
\end{discussionbox}

\begin{proposition}
\label{proposition:hiMiX}
When $X$ is as in Setup~\ref{setup:convexpolyomino}.
Write $h(t ) = 1 + h_1 t + h_2 t^2 + \cdots$ for the $h$-polynomial of
$\Bbbk[X]$.
Then $h_i = |\maxChains_i(X)|$.
\end{proposition}

\begin{proof}
Use~\cite[Theorems~4.1 and~1.1]{BjornerGarsiaStanley1982} with standard
grading (i.e. setting $t_i = t$ for all $i$) to see that the $h$-polynomial
of the Stanley Reisner ring of $\Delta(V(X)))$ is
\[
\sum_{i \in \naturals } 
|\maxChains_i(X)|t^i.
\]
The proposition now follows from Remark~\ref{remark:hPolyHibiRingSRRing}.
\end{proof}

\begin{discussionbox}
\label{discussionbox:JI}
Let $X$ be as in Setup~\ref{setup:convexpolyomino}.
Left-boundary vertices and bottom-boundary vertices are join-irreducible.
Let $p \in V(X )$; assume that $p$ is not a 
left-boundary vertex or a bottom-boundary vertex.
If $p \not \in \partial X$ then it is the 
top-right vertex of a cell in $X$, and hence is not join-irreducible.
If $p \in \partial X$ then $p$ is the bottom-left vertex of 
the unique cell containing it (i.e., the bottom element $\hat 0$ of $V(X)$) 
or the top-right vertex of the unique cell containing it
(i.e., the top element $\hat 1$ of $V(X)$); hence $p \not \in \JI(X )$.
Thus we have established that $\JI(X)$ is the union of the set of the
left-boundary vertices and of the set of the bottom-boundary vertices. 
The sets of the left-boundary vertices and of the bottom-boundary vertices
are totally ordered in $V(X)$. Therefore if $(p,p' )$ is a pair of
incomparable elements of $\JI(X )$, then one of them is a
left-boundary vertex and the other is a bottom-boundary vertex.
\end{discussionbox}

\section{Proof of the theorem}
\label{section:proof}

\begin{proposition}
\label{proposition:succDesc}
Let $\mu \in \maxChains(X )$ and $i \in \Des(\mu )$.
Write $\mu$ as a chain of order ideals 
$\hat 0 = I_0 \subsetneq I_1 \subsetneq \cdots \subsetneq I_{m+n} = \hat 1$
and $|I_i \minus I_{i-1}| = \{p_i\}$ with $p_i \in \JI(X )$.
Then
\begin{enumerate}

\item
\label{proposition:succDesc:incomp}
$p_i$ and $p_{i+1 }$ are incomparable;

\item
\label{proposition:succDesc:noSucc}
$i+1 \not \in \Des(\mu )$.
\end{enumerate}
\end{proposition}

\begin{proof}
\eqref{proposition:succDesc:incomp}:
Assume, by way of contradiction, that they are comparable. Then 
$p_i < p_{ i+1}$. Hence $\omega(p_i ) < \omega(p_{i+1 } )$, contradicting
the hypothesis that $i \in \Des(\mu )$.

\eqref{proposition:succDesc:noSucc}:
By way of contradiction, assume that $i+1 \in \Des(\mu )$.
Then, by~\eqref{proposition:succDesc:incomp}, $p_{i+1}$ and $p_{i+2 }$ are
incomparable. We see from Discussion~\ref{discussionbox:JI} and the
definition of the $p_i$ that $p_i < p_{i+2 }$. 
Therefore $\omega(p_i ) < \omega(p_{i+2} )$ contradicting the hypothesis
that $\omega(p_i ) > \omega(p_{i+1 } ) > \omega(p_{i+2} )$.
\end{proof}

\begin{proposition}
\label{proposition:psiInjective}
The function $\psi$ of~\eqref{equation:MXtoPowCX} is injective.
\end{proposition}

\begin{proof}
Let $\mu,\nu \in \maxChains(X )$ be such that $\psi(\mu) = \psi(\nu)$.
As earlier, write $\mu$ and $\nu$ as chains of order ideals of $\JI(X)$:
\begin{align*}
\mu & :  \hat 0 = 
I_0 \subsetneq I_1 \subsetneq \cdots \subsetneq I_{m+n} = \hat 1;
\\
\nu & :  \hat 0 = 
I'_0 \subsetneq I'_1 \subsetneq \cdots \subsetneq I'_{m+n} = \hat 1.
\end{align*}
For $1 \leq i \leq m+n$, write $I_i \minus I_{i-1} = \{p_i\}$ 
and $I'_i \minus I'_{i-1} = \{p'_i\}$ 
with $p_i , p'_i\in \JI(X )$.
We will prove by induction on $i$ that $I_i = I'_i$ for all 
$0 \leq i \leq m+n$. Since $I_0 = I'_0$, we may assume that $i>0$ and that
$I_j = I'_j$ for all $j<i$.

Assume, by way of contradiction, that $I_i \neq I'_i$. Then $I_{i-1}$
(which equals $I'_{i-1 }$) is the bottom-left vertex of a cell $C$. 
Without loss of generality, we may assume that $I_i$ is the top-left vertex
of $C$ and that $I'_{i }$ is the bottom-right vertex of $C$.
(In other words, $\mu$ goes up and $\nu$ goes to the right from $I_{i-1}$,
or equivalently, $p_i$ is a left-boundary vertex and $p'_i$ is a
bottom-boundary vertex.)

Let 
\begin{align*}
i_1 & = \min \{j > i : p'_i \in I_j \} - 1;
\\
i_2 & = \min \{j > i : p_i \in I'_j \} - 1.
\end{align*}
Then the edge 
$(I_{i_1-1}, I_{i_1})$ is vertical while $(I_{i_1}, I_{i_1+1})$ is
horizontal; this is the first time $\mu$ turns horizontal after $I_{i-1}$.
Let $C_1$ be the cell with 
$I_{i_1-1}$, $I_{i_1}$ and $I_{i_1+1}$ as the bottom-left, the top-left and
the top-right vertices respectively.
Similarly the edge 
$(I'_{i_2-1}, I'_{i_2})$ is vertical while $(I'_{i_2}, I'_{i_2+1})$ is
horizontal; this is the first time $\nu$ turns vertical after $I'_{i-1}$.
Let $C_2$ be the cell with 
$I'_{i_2-1}$, $I'_{i_2}$ and $I'_{i_2+1}$ as the bottom-left, the
bottom-right and the top-right vertices respectively.
(The possibility that $C_1 = C$ or $C_2=C$ has not been ruled out.)
See Figure~\ref{figure:munu} for a schematic showing 
the cells $C$, $C_1$ and $C_2$ and the chains $\mu$ and $\nu$.

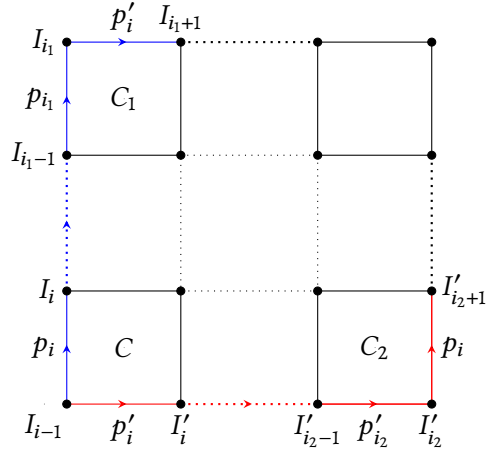
\begin{figure}
	\centering	
	\begin{tikzpicture}[scale=3]
	\draw[]  (.5,1.6)--(.5,1.1)--(0,1.1)
	(1.1,1.6)--(1.6,1.6)
	(1.6,.5)--(1.1,.5)--(1.1,0)
	(1.6,1.1)--(1.6,1.6)
	(0.5,0)--(0.5,.5)--(0,.5)
	(1.1,1.1)--(1.6,1.1)
	(1.1,1.1)--(1.1,1.6) ;
	
	\draw[blue] (0,0) -- (0,.5) node[ sloped, pos=0.5,  allow upside down]{\arrowIn};  
	\draw[thick, dotted, blue] (0,.5) -- (0,1.1) node[ sloped, pos=0.5,  allow upside down]{\arrowIn};  	
	\draw[blue] (0,1.1)--(0,1.6) node[ sloped, pos=0.5,  allow upside down]{\arrowIn};
	\draw[blue] (0,1.6)--(.5,1.6) node[ sloped, pos=0.5,  allow upside down]{\arrowIn};
	
	\draw[red] (0,0)--(.5,0) node[ sloped, pos=0.5,  allow upside down]{\arrowIn};  
	\draw[thick, dotted, red] (.5,0)--(1.1,0) node[ sloped, pos=0.5,  allow upside down]{\arrowIn};  	
	\draw[red] (1.1,0)--(1.6,0) node[ sloped, pos=0.5,  allow upside down]{\arrowIn};
	\draw[red] (1.6,0)--(1.6,.5) node[ sloped, pos=0.5,  allow upside down]{\arrowIn};	
	
	\draw[red]   (1.1,0)--(1.6,0)--(1.6,.5);
	
	\draw[thick, dotted ]  (.5,1.6)--(1.1,1.6) (1.6,.5)--(1.6,1.1) ; 
	
	\draw[dotted]   (1.1,.5)--(.5,.5)--(.5,1.1)
	(1.1,.5)--(1.1,1.1) (.5,1.1)--(1.1,1.1)  ;
	
	\filldraw[black] (0,0) circle (.5pt) node[anchor=north]  {}; 
	\filldraw[black] (.5,0) circle (.5pt) node[anchor=north] {$I_i'$}; 
	\filldraw[black] (0,.5) circle (.5pt) node[anchor=east] {$I_i$}; 
	\filldraw[black] (.5,.5) circle (.5pt) node[anchor=west] {};  
	\filldraw[black] (.5,1.1) circle (.5pt) node[anchor=north] {}; 
	\filldraw[black] (1.1,.5) circle (.5pt) node[anchor=north] {}; 
	\filldraw[black] (0,1.1) circle (.5pt) node[anchor=east] {$I_{{i_1}-1}$};  
	
	\filldraw[black] (1.6,0.5) circle (.5pt) node[anchor=west] {$I'_{{i_2}+1}$}; 
	\filldraw[black] (0.5,1.6) circle (.5pt) node[anchor=south] {$I_{{i_1}+1}$}; 
	\filldraw[black] (1.6,1.6) circle (.5pt) node[anchor=west] {};  
	
	\filldraw[black] (1.6,0) circle (.5pt) node[anchor=north] {$I'_{{i_2}}$}; 
	
	\filldraw[black] (-.1,0) circle (0pt) node[anchor=north] {$I_{i-1}$};	
	
	\filldraw[black] (0,1.6) circle (.5pt) node[anchor=east] {$I_{{i_1}}$};
	\filldraw[black] (1.1,1.1) circle (.5pt) node[anchor=east] {}; 	
	\filldraw[black] (1.1,0) circle (.5pt) node[anchor=north] {$I'_{{i_2}-1}$}; 
	\filldraw[black] (1.1,1.6) circle (.5pt) node[anchor=west] {}; 
	\filldraw[black] (1.6,1.1) circle (.5pt) node[anchor=east] {}; 
	
	\filldraw[black] (.25,0) circle (0pt) node[anchor=north] {$p_i'$}; 
	\filldraw[black] (.25,1.6) circle (0pt) node[anchor=south] {$p_i'$}; 
	\filldraw[black] (0,0.25) circle (0pt) node[anchor=east] {$p_i$}; 
	\filldraw[black] (1.6,0.25) circle (0pt) node[anchor=west] {$p_i$}; 
	\filldraw[black] (0,1.35) circle (0pt) node[anchor=east] {$p_{i_1}$};          		
	\filldraw[black] (1.35,0) circle (0pt) node[anchor=north] {$p'_{i_2}$};   
	
	\filldraw[] (.25,0.25) circle (0pt) node[anchor=center]  {$C$};	
	\filldraw[] (1.35,.25) circle (0pt) node[anchor=center]  {$C_2$};	
	\filldraw[] (.25,1.35) circle (0pt) node[anchor=center]  {$C_1$};		
	
	\end{tikzpicture}
\caption{$C$, $C_1$, $C_2$, $\mu$ (blue) and $\nu$ (red) from the proof of 
Proposition~\ref{proposition:psiInjective}.}
\label{figure:munu}
\end{figure}

We now prove a sequence of statements from which the proposition follows.

\begin{asparaenum}

\item
\label{proposition:psiInjective:COneNotInPsiMu}
If $C_1 \not \in  \psi(\mu )$, then $C_2 \in \psi(\nu )$. 
Proof: 
Note that $p_{i_1+1} = p'_{i}$ and $p'_{i_2+1 } = p_i$.
Since $C_1 \not \in  \psi(\mu )$,
we see that 
\[
\omega(p'_{i}) = \omega(p_{i_1+1}) > \omega(p_{i_1}) \geq 
\omega(p_{i}),
\]
where the last inequality follows from noting that 
$p_i < \cdots < p_{i_1 }$ since they are left-boundary vertices.
Therefore, in the chain $\nu$, we have
\[
\omega(p'_{i_2}) \geq \omega(p'_{i}) > 
\omega(p_{i})= \omega(p'_{i_2+1}), 
\]
i.e., $i_2 \in \Des(\nu )$. Hence $C_2 \in \psi(\nu )$.

\item
\label{proposition:psiInjective:CTwoNotInPsiMu}
If $C_2 \not \in  \psi(\nu )$, then $C_1 \in \psi(\mu )$. 
Immediate from \eqref{proposition:psiInjective:COneNotInPsiMu}.

\item
\label{proposition:psiInjective:mu}
If $C_1 \neq C$ then $C \not \in \psi(\mu)$ and $C_1 \not \in \psi(\nu)$.
Proof:
Note that $\mu$ does not pass through the top-right vertex of $C$ 
 and that $\nu$ does not pass through the bottom-left vertex of $C_1$.

\item
\label{proposition:psiInjective:nu}
If $C_2 \neq C$ then $C \not \in \psi(\nu)$ and $C_2 \not \in \psi(\mu)$.
Proof:
Note that $\nu$ does not pass through the top-right vertex of $C$ 
 and that $\mu$ does not pass through the bottom-left vertex of $C_1$.

\item
\label{proposition:psiInjective:COneDiffC}
If $C_1 \neq C$, then $\psi(\mu ) \neq \psi(\nu )$. 
Proof:
If $C_1 \in \psi(\mu )$, use~\eqref{proposition:psiInjective:mu} to see
that
\[
C_1 \in \psi(\mu ) \minus \psi(\nu ).
\]
Now assume that $C_1 \not \in \psi(\mu )$. 
Then $C_2 \in \psi(\nu )$
by~\eqref{proposition:psiInjective:COneNotInPsiMu}. 
If $C_2 = C$, then $C_2 \not \in \psi(\mu )$
by~\eqref{proposition:psiInjective:mu};
otherwise, $C_2 \not \in \psi(\mu )$
by~\eqref{proposition:psiInjective:nu}.

\item
\label{proposition:psiInjective:CTwoDiffC}
If $C_2 \neq C$, then $\psi(\mu ) \neq \psi(\nu )$. 
Proof:
If $C_2 \in \psi(\nu )$, use~\eqref{proposition:psiInjective:nu} to see
that
\[
C_2 \in \psi(\nu ) \minus \psi(\mu ).
\]
Now assume that $C_2 \not \in \psi(\nu )$. 
Then $C_1 \in \psi(\mu )$
by~\eqref{proposition:psiInjective:CTwoNotInPsiMu}. 
If $C_1 = C$, then $C_1 \not \in \psi(\nu )$
by~\eqref{proposition:psiInjective:nu};
otherwise, $C_1 \not \in \psi(\nu )$
by~\eqref{proposition:psiInjective:mu}.

\item
\label{proposition:psiInjective:atmostone}
$C$ belongs to at most one of $\psi(\mu )$ and $\psi(\nu )$.
Proof:
Suppose $C \in \psi(\mu)$. Then $i_1 = i+1$, $p_{i_1 } = p'_i$
and $\omega(p_i ) > \omega(p'_i )$. For $C$ to belong to $\psi(\nu )$, we
need that $I'_{i+1} = I_{i+1}$ (i.e., $\mu$ and $\nu$ are the same up to
$i+1$, except at $i$); for this to hold, it is necessary that $p'_{i+1 } =
p_i$, but then $i \not \in \Des(\nu)$. The other case is proved similarly.

\item
\label{proposition:psiInjective:BothEqC}
If $C_1 = C_2  = C$ then $\psi(\mu ) \neq \psi(\nu )$. 
Proof:
By~\eqref{proposition:psiInjective:atmostone}, it suffices to show that 
$C \in \psi(\mu )$ or $C \in \psi(\nu )$. This follows
from~\eqref{proposition:psiInjective:COneNotInPsiMu}
and~\eqref{proposition:psiInjective:CTwoNotInPsiMu}.
\end{asparaenum}

The proposition is proved by~\eqref{proposition:psiInjective:COneDiffC}, 
~\eqref{proposition:psiInjective:CTwoDiffC},
and~\eqref{proposition:psiInjective:BothEqC}.
\end{proof}

\begin{proposition}
\label{proposition:rookConfig}
Let $k \in \naturals$ and $\mu \in \maxChains_k(X )$.
Then $\psi(\mu )$ is a $k$-rook configuration in $X$.
\end{proposition}

\begin{proof}
Since $|\psi(\mu ) | = k$, it suffices to note that the cells of 
$\psi(\mu )$ are in distinct rows and columns. This follows from
Proposition~\ref{proposition:succDesc}\eqref{proposition:succDesc:noSucc}.
\end{proof}

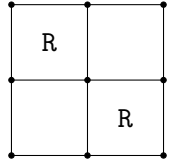
\begin{figure}
	\centering	
	\begin{tikzpicture}[scale=2]
	\draw[] (0,0)--(0,1)--(1,1)--(1,0)--(0,0) (0,.5)--(1,.5) (.5,0)--(.5,1);
	
	\filldraw[black] (0,0) circle (.5pt) node[anchor=north]  {}; 
	\filldraw[black] (.5,0) circle (.5pt) node[anchor=north] {}; 
	\filldraw[black] (0,.5) circle (.5pt) node[anchor=east] {}; 
	\filldraw[black] (.5,.5) circle (.5pt) node[anchor=west] {};  
	\filldraw[black] (1,1) circle (.5pt) node[anchor=north]  {}; 
	\filldraw[black] (.5,1) circle (.5pt) node[anchor=north] {}; 
	\filldraw[black] (1,.5) circle (.5pt) node[anchor=east] {}; 
	\filldraw[black] (0,1) circle (.5pt) node[anchor=west] {};  
	\filldraw[black] (1,0) circle (.5pt) node[anchor=east] {}; 
	
	\filldraw[] (.75,0.25) circle (0pt) node[anchor=center]  {\texttt{R}};		
	\filldraw[] (.25,0.75) circle (0pt) node[anchor=center]  {\texttt{R}};		
	\end{tikzpicture}
	\caption{$2$-rook (denoted by \texttt{R})
    configuration in a non-thin polyomino.}
    \label{figure:twoRooks}
\end{figure}

\begin{proof}[Proof of Theorem~\protect{\ref{theorem:RRconj}}]
For each $i \in \naturals$, $h_i = |\maxChains_i(X)|$
by Proposition~\ref{proposition:hiMiX}.
By Propositions~\ref{proposition:psiInjective}
and~\ref{proposition:rookConfig} we see that $h_i \leq r_i$ for all $i$.
Since $X$ is not thin, $X$ contains a $2$-rook configuration as in
Figure~\ref{figure:twoRooks}.
Such a rook configuration cannot be in the image of $\psi$. 
Hence $h_2 < r_2$.
\end{proof}

Using results of~\cite{EneHerzogQureshiRomeoLconvex2021}, 
we can extend our result to $L$-convex polyominoes as follows.
Let $X$ be an $L$-convex polyomino. Then there exists a polyomino $X^*$
(the Ferrer diagram projected by $X$, in the sense 
of~\cite{EneHerzogQureshiRomeoLconvex2021}) such that 
\begin{enumerate}

\item
$X^*$ is a convex polyomino such that $V(X^*)$ is a sublattice of
$\naturals^2$ (since $X^*$ is a Ferrer diagram);

\item
If $X$ is not thin, then $X^*$ is not thin;

\item
$X$ and $X^*$ have the same rook
polynomial~\cite[Lemma~2.4]{EneHerzogQureshiRomeoLconvex2021};

\item
$\Bbbk[X]$ and $\Bbbk[X^*]$ are isomorphic to each 
other~\cite[Theorem~3.1]{EneHerzogQureshiRomeoLconvex2021}, so they have the
same $h$-polynomial.
\end{enumerate}

Thus we get:

\begin{corollary}
\label{corollary:Lconvex}
Let $X$ be an $L$-convex polyomino that is not thin.
Let $h(t) = 1 + h_1 t + h_2 t^2 + \cdots$ be the $h$-polynomial of 
$\Bbbk[X]$ and $r(t) = 1 + r_1 t + r_2 t^2 + \cdots$ be the rook polynomial
of $X$.
Then $h_2 < r_2$.
\end{corollary}

\def\cfudot#1{\ifmmode\setbox7\hbox{$\accent"5E#1$}\else
  \setbox7\hbox{\accent"5E#1}\penalty 10000\relax\fi\raise 1\ht7
  \hbox{\raise.1ex\hbox to 1\wd7{\hss.\hss}}\penalty 10000 \hskip-1\wd7\penalty
  10000\box7}

\end{document}